\documentclass[12pt,a4paper]{article}

\usepackage{amsmath,amsfonts,amssymb,amsthm,graphicx,authblk,setspace}     
\usepackage{color}
\usepackage{booktabs}

\newcommand{\bs}{\bigskip}

\newcommand{\calA}{\mathcal{A}}

\newcommand{\calF}{\mathcal{F}}

\newcommand{\al}{\alpha}

\newcommand{\tha}{\theta}

\newcommand{\Up}{\Upsilon}

\newcommand{\ze}{\zeta}

\newcommand{\E}{\mathbb{E}}
\newcommand{\Ex}{\mathbb{E}}

\newcommand{\N}{\mathbb{N}}
\newcommand{\Pb}{\mathbb{P}}
\newcommand{\R}{\mathbb{R}}

\newcommand{\Z}{\mathbb{Z}}

\newcommand{\bfi}{\mathbf{i}}
\newcommand{\bfj}{\mathbf{j}}

\newcommand{\bfu}{\mathbf{u}}

\newcommand{\half}{{\textstyle \frac{1}{2}}}

\newcommand{\eqdist}{\buildrel{{\rm d}}\over =}

\newtheorem{thm}{Theorem}
\newtheorem{assp}{Assumption}
\newtheorem{lem}{Lemma}
\newtheorem{cor}{Corollary}
\newtheorem{definition}{Definition}
\newtheorem{rema}{Remark}

\begin{document}

\title{\normalsize \textbf{MODULUS OF CONTINUITY OF A CLASS OF MONOFRACTAL PROCESSES}}

\author{ \ \normalsize Geoffrey Decrouez$^1$, Ben Hambly$^2$\ and\ Owen Dafydd Jones$^1$} 
  
\footnotetext[1]{Department of Mathematics and Statistics, The University of Melbourne, VIC 3010, Australia}

\footnotetext[2]{Mathematical Institute, University of Oxford, Oxford, UK}  

\date{}
\maketitle

\medskip\noindent {\normalsize \textbf{Abstract.}}
We derive the modulus of continuity of a class of processes called Canonical Embedded Branching Processes (CEBP), recently introduced by Decrouez and Jones \cite{DecrJ??}, and we establish their monofractal character. CEBP provide a rich class of processes, including the Brownian motion as a particular case. The techniques developed in this study follow the steps of Barlow and Perkins on Brownian motion on a Sierpinski gasket \cite{BarlP88}, though complications arise here since CEBP are not Markovian in general.


\bs
\begin{center}{\bf 1.  INTRODUCTION}\end{center}

The local fluctuations of a process $X$ can be described using the local H\"older exponent $h_X(t)$, defined as
\[
h_X(t):= \liminf_{\epsilon\rightarrow 0} \frac{1}{\log \epsilon}\log \sup_{|u-t|<\epsilon} |X(u)-X(t)|.
\]
When $h_X(t)$ is constant all along the sample path with probability $1$, $X$ is said to be monofractal.
In contrast, there exist processes whose H\"older exponent behaves erratically, whereby in any interval of positive length we find a range of different exponents.
For such processes, it is in practice impossible to estimate $h_X(t)$ for all $t$, due to the finite precision of the data.
Instead, we use the multifractal or Hausdorff spectrum $D_X$, a global description of its local fluctuations.
$D_X(h)$ is defined as the Hausdorff dimension of the set of points with a given H\"older exponent $h$.
For monofractal processes, $D_X(h)$ degenerates to a single point at some $h=H$ (so $D_X(H)=1$, and the convention is to set $D_X(h)=-\infty$ for $h\neq H$).
When the spectrum is non trivial for a range of values of $h$, the process is said to be multifractal.

Recently Decrouez and Jones \cite{DecrJ??} described a new class of processes, called Canonical Embedded Branching Process (CEBP) processes.
CEBP are defined using the crossing tree, an ad-hoc space-time description of the process, and are such that the spatial component of their crossing tree is a Galton-Watson branching process.
The contribution of this study is to obtain the modulus of continuity for CEBP processes, and hence show that they are monofractal.

\newpage
\begin{center}{\bf 2.  THE CANONICAL EMBEDDED BRANCHING PROCESS (CEBP)}\end{center}
\label{sec:EBP}

Let $X: \R^+\rightarrow\R$ be a continuous process, with $X(0)=0$.
For $n\in\Z$ we define level $n$ passage times $T_{k}^{n}$ by putting $T_0^n = 0$ and
\begin{equation*}
T_{k+1}^{n} = \inf\{t>T_{k}^{n}~|~X(t)\in  2^{n}\Z,~X(t)\not = X(T_{k}^{n}) \}.
\end{equation*}
The $k$-th level $n$ (equivalently scale $2^n$) crossing $C_k^n$ is the sample path from $T_{k-1}^{n}$ to $T_{k}^{n}$.
That is, $C_k^n = \{ (t,X(t))\mid T_{k-1}^n\leq t < T_k^n\}$.

\begin{figure}
\centering
\includegraphics[width=12cm]{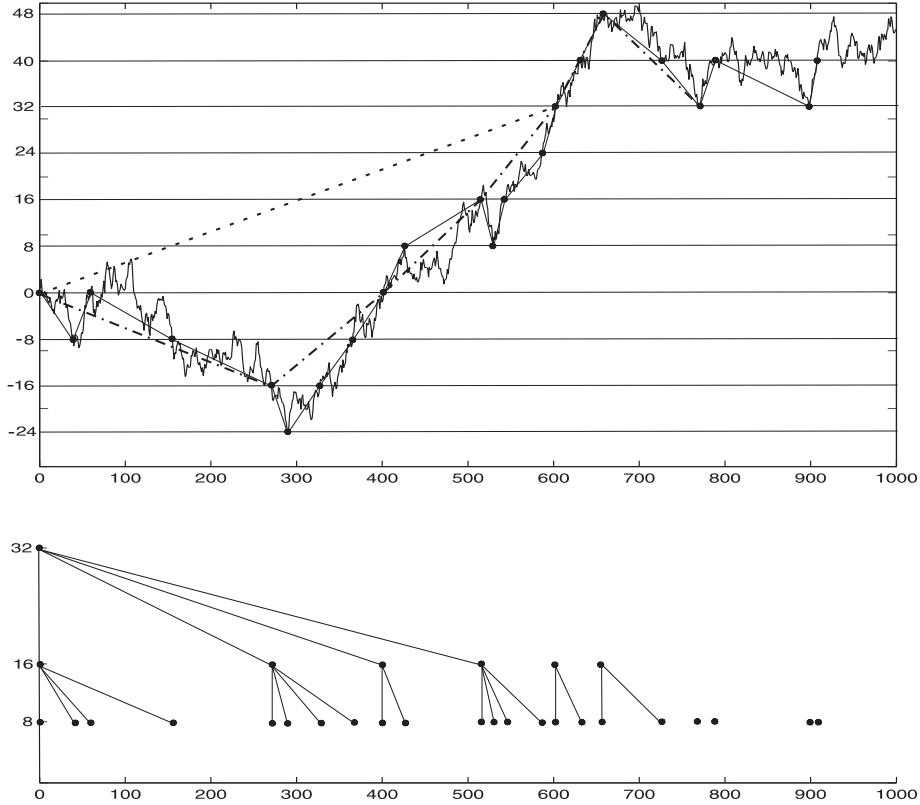}
\caption{A section of sample path and levels 3, 4 and 5 of its crossing tree.
In the top frame we have joined the points $T^n_k$ at each level, and in the bottom frame we have identified the $k$-th level $n$ crossing with the point $(2^n, T^n_{k-1})$ and linked each crossing to its subcrossings.}%
\label{crossingtree}
\end{figure}

When passing from a coarse scale to a finer one, we decompose each level $n$ crossing into a sequence of level $n-1$ crossings.
To define the crossing tree, we associate nodes with crossings, and the children of a node are its subcrossings.
The crossing tree is illustrated Figure \ref{crossingtree}, where the level 3, 4 and 5 crossings of a given sample path are shown.

In addition to indexing crossings be their level and position within each level, we will also use a tree indexing scheme.
Let $\emptyset$ be the root of the tree, representing the first level 0 crossing.
The first generation of children (which are level $-1$ crossings, of size $1/2$) are labelled by $i$, $1\leq i \leq Z_\emptyset$, where $Z_\emptyset$ is the number of children of $\emptyset$.
The second generation (which are level $-2$ crossings, of size $1/4$) are then labelled $ij$, $1\leq j\leq Z_i$, where $Z_i$ is the number of children of $i$.
More generally, a node is an element of $U = \cup_{m\geq 0}\N^{m}$ and a branch is a couple $(\bfu,\bfu j)$ where $\bfu\in U$ and $j\in\N$.
The length of a node $\bfi=i_1\hdots i_m$ is $|\bfi|=m$, and the $k$-th element is $\bfi[k] = i_k$.
If $|\bfi|>m$, $\bfi|_m$ is the curtailment of $\bfi$ after $m$ terms.
Conventionally $|\emptyset|=0$ and $\bfi|_0=\emptyset$.
A tree $\Up$ is a set of nodes, that is a subset of $U$, such that
\begin{itemize}
\item $\emptyset\in \Up$
\item If a node $\bfi$ belongs to the tree then every ancestor node $\bfi|_k$, $k\leq |\bfi|$, belongs to the tree
\item If $\bfu\in \Up$, then $\bfu j\in \Up$ for $j = 1, \ldots, Z_\bfu$ and $\bfu j \not\in \Up$ for $j > Z_\bfu$, where $Z_\bfu$ is the number of children of $\bfu$.
\end{itemize}

Let $\Up_m$ be the $m$-th generation of the tree, that is the set of nodes of length $m$.
(These are level $-m$ crossings, of size $2^{-m}$.)
Define $\Up_\bfi= \{\bfj\in\Up \,|\, |\bfj|\geq |\bfi| \textrm{ and } \bfj|_{|\bfi|}=\bfi\}$.
The boundary of the tree is given by $\partial \Up = \{\bfi \in \N^\N \,|\, \forall m\geq 0, \bfi|_m\in \Up \}$.
Let $\psi(\bfi)$ be the position of node $\bfi$ within generation $|\bfi|$, so that crossing $\bfi$ is just $C^{-|\bfi|}_{\psi(\bfi)}$.
The nodes to the left and right of $\bfi$, corresponding to the crossings $C^{-|\bfi|}_{\psi(\bfi)-1}$ and $C^{-|\bfi|}_{\psi(\bfi)+1}$, will be denoted $\bfi-$ and $\bfi+$.
In general when we have quantities associated with crossings we will use tree indexing and level/position indexing interchangeably.
So $Z_\bfi = Z^{-|\bfi|}_{\psi(\bfi)}$, $T_\bfi = T^{-|\bfi|}_{\psi(\bfi)}$, etc.

Let $\al^n_k \in \{+, -\}$ be the orientation of $C^n_k$, $+$ for up and $-$ for down, and let $A^n_k$ be the vector given by the orientations of the subcrossings of $C^n_k$.
Let $D^n_k = T^n_k - T^n_{k-1}$ be the duration of $C^n_k$.
Clearly, to reconstruct the process we only need $\al^n_k$ and $D^n_k$ for all $n$ and $k$.
The $\al^n_k$ encode the spatial behaviour of the process, and the $D^n_k$ the temporal behaviour.
Our definition of an EBP is concerned with the spatial component only.

\begin{definition}
A continuous process $X$ with $X(0) = 0$ is called an Embedded Branching Process (EBP) process if for any fixed $n\in\Z$, conditioned on the crossing orientations $\al^n_k$ the random variables $A^n_k$ are all mutually independent, and $A^n_k$ is conditionally independent of all $A^m_j$ for $m > n$.
In addition we require that $\{A^n_k \,|\, \al^n_k = i\}$ are identically distributed, for $i = +, -$.

That is, an EBP process is such that if we take any given crossing, then count the orientations of its subcrossings at successively finer scales, we get a (supercritical) two-type Galton-Watson process, where the types correspond to the orientations.
\end{definition}

Subcrossing orientations have a particular structure.
A level $n$ up crossing is from $k2^n$ to $(k+1)2^n$, a down crossing is from $k2^n$ to $(k-1)2^n$.
The level $n-1$ subcrossings that make up a level $n$ parent crossing consist of {\em excursions} (up-down and down-up pairs) followed by a {\em direct crossing} (down-down or up-up pairs), whose direction depends on the parent crossing: if the parent crossing is up, then the subcrossings end up-up, otherwise, they end down-down.
Let $Z^n_k$ be the length of $A^n_k$, that is, the number of subcrossings of $C^n_k$.
The number of up and down subcrossings will be written $Z^{n+}_k$ and $Z^{n-}_k$ respectively.
Clearly, each of the $Z_k^n-2$ first entries of $A_k^n$ come in pairs, each pair being up-down or down-up. The last two components are either the pair up-up or down-down, depending on $\alpha_k^n$.
Thus, given $\al^n_k = +$, we must have $Z^{n+}_k = \half Z^n_k + 1$ and $Z^{n-}_k = \half Z^n_k - 1$, and conversely given $\al^n_k = -$.

Let $\calA$ be the space of possible orientations.
That is, $a \in \calA$ consists of some number of pairs, $+-$ or $-+$, then a single pair $++$ or $--$.
Given an EBP process, for the offspring type distributions we write $p^+_{A}(a) = \Pb(A^n_k = a \,|\, \al^n_k = +)$ and $p^-_{A}(a) = \Pb(A^n_k = a \,|\, \al^n_k = -)$, for $a \in \calA$.
Let $\mu^+ = \E(Z^n_k\,|\,\al^n_k=+)$, $\mu^- = \E(Z^n_k \,|\, \al^n_k=-)$ and $\mu = \half(\mu^+ + \mu^-)$, then the mean offspring matrix is given by
\[
M := \E \left(
\begin{array}{cc}
(Z^{n+}_k|\al^n_k = +) & (Z^{n-}_k|\al^n_k = +) \\
(Z^{n+}_k|\al^n_k = -) & (Z^{n-}_k|\al^n_k = -)
\end{array}
\right)
= \left(
\begin{array}{cc}
\half \mu^+ + 1 & \half \mu^+ - 1 \\
\half \mu^- - 1 & \half \mu^- + 1
\end{array}
\right)
\]

To proceed we need to make some assumptions about $p^\pm_{A}$.

\begin{assp}\label{AssGW}
$\mu^+, \mu^- > 2$ and $\E(Z^{ni}_k \log Z^{ni}_k \,|\, \al^n_k = j) < \infty$ for $i, j = \pm$
\end{assp}

The first of these assumptions ensures that $M$ is strictly positive with dominant eigenvalue $\mu > 2$, and corresponding left eigenvector $(\half, \half)$.
The corresponding right eigenvector is $((\mu^+-2)/(\mu-2), (\mu^--2)/(\mu-2))^T$.
The second assumption is the usual condition for the normed limit of a supercritical Galton-Watson process to be non-trivial.

\begin{thm}\label{CEBP} (\cite{DecrJ??} Theorems 2.1 and 2.2)
For any offspring orientation distributions $p^\pm_{A}$ satisfying Assumption \ref{AssGW}, there exists a corresponding continuous EBP process $X$ defined on $\R_+$.
Moreover, we can choose $X$ so that for each $n\in\Z$, conditioned on the crossing orientations $\al^n_k$, the crossing durations $D^n_k$ are all mutually independent, and $D^n_k$ is conditionally independent of all $A^m_j$ for $m > n$.
Also, $\E(D^n_k \,|\, \al^n_k = \pm) = \mu^{n} (\mu^\pm - 2)/(\mu - 2)$, and the distribution of $\mu^{-n}D^n_k$ depends only on $\al^n_k$.

Up to finite dimensional distributions, $X$ is the unique such EBP with offspring orientation distributions $p^\pm_{A}$.
That is, for any other EBP process $Y$ with offspring orientation distributions $p^\pm_A$ and crossing durations as above, we have $(X(t_1), \ldots, X(t_k)) \eqdist (Y(t_1), \ldots, Y(t_k))$ for any $0 \leq t_1 < t_2 < \cdots < t_k$.
Accordingly we call $X$ the Canonical EBP (CEBP) process with these offspring distributions.

The two-type Galton-Watson process defined by the orientation distributions $p^\pm_A$ is supercritical, and when scaled by its mean converges to $(\half, \half)W^i$, where $i = \pm$ is the type of the first individual.
For the Canonical EBP we have that the crossing duration $D^n_k \eqdist \mu^{n} W^i$, where $i = \al^n_k$.

We also observe that $X$ is discrete scale-invariant: let $H=\log2/\log\mu$, then for all $c\in\{\mu^n, n\in\Z\}$,
\begin{equation}\label{eq:holder}
X(t) \stackrel{fdd}{=} c^{-H} X(ct),
\end{equation}
where $\stackrel{fdd}{=}$ denotes equality for finite dimensional distributions.
$H = \log\mu/\log 2$ is known as the Hurst index.
\end{thm}

To simplify our lives somewhat, we will from here on restrict ourselves to EBP for which the embedded Galton-Watson process is single type.
That is, the number of subcrossings does not depend on the orientation of the parent crossing.
We will also suppose that, given the number of subcrossings, the excursions are independent and equally likely to be up-down and down-up.
That is, we will make the following assumption.

\begin{assp}\label{AssSym}
The distribution of $Z^n_k$ is independent of $\al^n_k$, and
\[
p^+_A(\cdots ++) = p^-_A(\cdots --) = 2^{-z}\Pb(Z^0_1 = 2(z + 1)),
\]
where $\cdots$ represents a combination of $z$ pairs, each either $+-$ or $-+$.
\end{assp}

Write $Z$ for a generic r.v.\ with the same distribution as the $Z^n_k$

When Assumptions \ref{AssGW} and \ref{AssSym} hold, the CEBP is completely specified by the subcrossing family size distribution $Z$.
Starting with any single crossing, if we just count the number of subcrossings at successively finer scales, then we get a single-type Galton-Watson process with offspring distribution $Z$.
Clearly $\mu^+ = \mu^- = \mu = \Ex Z$.
Let $W$ be the limit of the (supercritical) Galton-Watson process with offspring distribution $Z$, normed by its mean, then $\Ex W = 1$ and for any given $n$ the crossing durations $D^n_k$ are all distributed as $\mu^n W$.

\begin{center}{\bf 3.  THE MODULUS OF CONTINUITY OF A CEBP}\end{center}

The goal of this section is to establish that CEBP are monofractal processes with H\"older exponent $H=\log 2/\log\mu\in(0,1)$.
Our approach is based on that used by Barlow \& Perkins \cite{BarlP88} to obtain the modulus of continuity of Brownian motion on the Sierpinski gasket, though complications arise because CEBP are not in general Markovian.

The basic idea is to use bounds on the crossing durations $D^n_k$ to control how fast the process can move away from a given point.
We will take Assumptions \ref{AssGW} and \ref{AssSym} to hold throughout.
Let $X$ be the CEBP determined by the subcrossing number distribution $Z$.

We have the following from Biggins \& Bingham \cite{BiggB93}.

\begin{lem}\label{leftW}
Suppose that Assumptions \ref{AssGW} and \ref{AssSym} hold.
There exists strictly positive constants $c_1$, $c_2$, and $c_3$ such that for all $x>0$
$$
\exp\left( -c_1 x^{-H/(1-H)} \right)
\leq \Pb(W<x)
\leq c_2 \exp\left( -c_3 x^{-H/(1-H)} \right).
$$
\end{lem}

Let $T^n_k(s)$ be the $k+1$-st level $n$ crossing time greater than or equal to $s$, for $k \geq 0$.
So if $s$ is a level $n$ crossing time then $T^n_0(s) = 0$.
The previous lemma gave us a bound on the duration of a crossing.
The next lemma gives a lower bound on the time remaining in the current crossing.
To establish this result we will need to make a further modest restriction to the class of CEBP we consider.

\begin{assp}\label{AssZ}
We assume that the subcrossing number distribution $Z$ is such that there exists a $\ze$ such that for all $y$
\[
Z + \ze \geq_{st} Z - y \,|\, Z > y.
\]
Here $\geq_{st}$ denotes stochastic domination.
That is, for all $y$ and $z$,
\[
\Pb( Z - y > z \,|\, Z > y) \leq \Pb( Z + \ze > z).
\]
\end{assp}

This condition clearly holds for $Z$ bounded, and for $Z$ that are NBU (New Better than Used, in which case $\ze = 0$).
Examples of NBU distributions include the negative binomial with shape $\geq 1$ and the Poisson.

Let $\{ \calF_s \}_{s \geq 0}$ be the filtration generated by $X$.

\begin{lem}\label{leftT}
Suppose that Assumptions \ref{AssGW}, \ref{AssSym} and \ref{AssZ} hold, then there exists constants $c_4, c_5 >0$ such that for all $x>0$ and $n\in\Z$,
$$
\Pb( T^n_0(s) - s \leq x \,|\, \calF_s) \geq c_4 \exp\left( -c_5(\mu^{-n} x)^{-H/(1-H)} \right).
$$
\end{lem}

\begin{proof}
Note first that
\[
T^{n+1}_0(s) = T^n_0(s) + \sum_{i=1}^{Z^{n+1}(s)} \mu^{n} W(i),
\]
where $Z^{n+1}(s) \geq 0$ is the number of level $n$ crossings from $T^n_0(s)$ to $T^{n+1}_0(s)$, and the $W(i)$ are independent and distributed as $W$.
If $Z^{n+1}(s)$ is not zero then, conditioned on $\calF_s$, it will be distributed as $Z - y \,|\, Z > y$, where $Z$ has the subcrossing number distribution, and $y \geq 0$ is the number of level $n$ crossings from the current level $n+1$ crossing which have already happened by time $T^n_0(s)$, including the current level $n$ crossing. Notations are illustrated in Figure~\ref{lem27}.
\begin{figure}
\centering
\includegraphics[width=12cm]{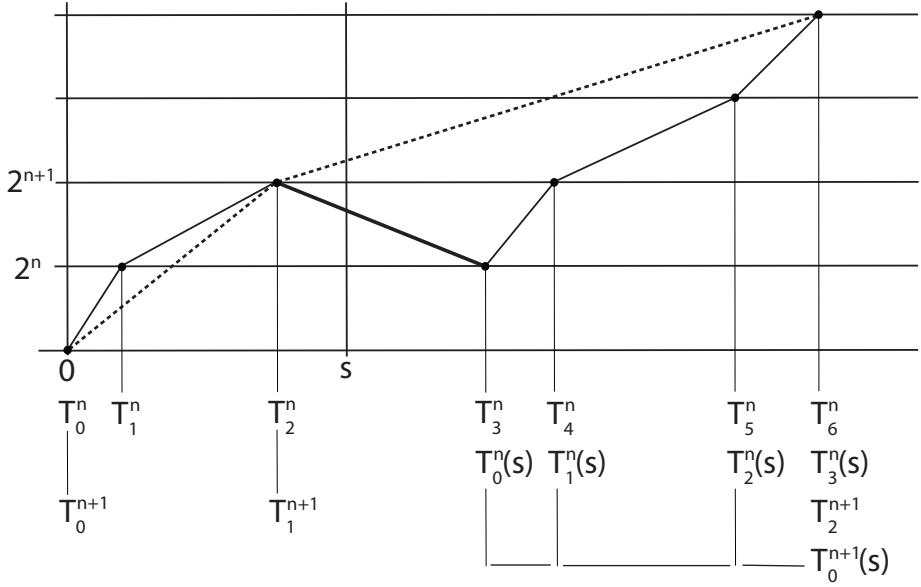}
\caption{Notations in Lemma~\ref{leftT}. The current level $n$ crossing at time $s$ is in bold. The first level $n$ crossing time greater than $s$, denoted $T_0^n(s)$, corresponds to $T_3^n$, and the first level $(n+1)$ crossing time greater than $s$ is $T_0^{n+1}(s)=T_6^n$. There are $Z^{n+1}(s)=3$ level $n$ crossings from $T_0^n(s)$ to $T_0^{n+1}(s)$, and $y=1$ level $n$ crossing from the current level $(n+1)$ crossing that has already happened at time $s$, including the current crossing in bold.}%
\label{lem27}
\end{figure}

Thus, from our assumption on $Z$, conditioning on $\calF_s$ we have
\begin{eqnarray*}
T^{n+1}_0(s) - T^n_0(s) &\leq_{st}& \sum_{i=1}^{Z + \ze} \mu^{n}W(i) \\
&\eqdist& \mu^{n+1}W(0) + \sum_{i=1}^\ze \mu^{n}W(i).
\end{eqnarray*}
As $n \downarrow -\infty$ we have $T^{n}_0(s) \downarrow s$ (this follows directly from \cite{O80} Theorem 1), whence
\[
T^n_0(s) - s \leq_{st} \mu^n W(0,n) + \sum_{k=-\infty}^{n-1} \sum_{i=0}^\ze \mu^{k} W(i,k),
\]
where the $W(i,k)$ are i.i.d.\ with distribution $W$.
Thus, for any $\tha > 0$ we have
\begin{equation}\label{lb1}
\Ex( e^{-\theta(T^n_0(s) - s)} \,|\, \calF_s)
\geq \Ex e^{-\tha \mu^n W} \prod_{k=-\infty}^{n-1} \left( \E e^{-\theta \mu^k W} \right)^\ze.
\end{equation}
Using the lower bound in Lemma \ref{leftW} for the left tail of $W$, it follows that for any $x>0$
\begin{align*}
\Ex e^{-\theta \mu^k W}
&\geq e^{-\theta x} \Pb(\mu^k W <x) \\
&\geq \exp\left( -\theta x - c_1(x\mu^{-k})^{-H/(1-H)} \right).
\end{align*}
For $x=(c_1/\theta)^{1-H}\mu^{kH}$ we get $\Ex e^{-\theta \mu^k W} \geq \exp\left( -c_2(\theta\mu^k)^H\right)$.
Plugging this into (\ref{lb1}) yields the bound
$$
\Ex( e^{-\theta(T^n_0(s) - s)} \,|\, \calF_s) \geq \exp\left( -c_3(\mu^n\theta)^H\right).
$$
Applying Lemma 4.1 in \cite{BarlP88} with $\theta = (c_4/x)^{1/(1-H)}2^{(2+k)/(1-H)}$, and readjusting the constants, gives the result.
\end{proof}

\begin{lem}\label{cebp}
Suppose that Assumptions \ref{AssGW}, \ref{AssSym} and \ref{AssZ} hold, then there exist constants $c_6, \ldots, c_9>0$ such that for all $\lambda>0$ and any $s, t \geq 0$,
\begin{eqnarray*}
\lefteqn{ c_6 \exp\left(-c_7 (\lambda^{1/H}/t)^{H/(1-H)}\right) \leq \Pb(|X(s+t)-X(s)|>\lambda\mid\calF_s)} \\
&& \leq \Pb\left(\sup_{0\leq u\leq t}
|X(s+u)-X(s)|>\lambda\mid\calF_s\right) \leq c_8 \exp\left(-c_9 (\lambda^{1/H}/t)^{H/(1-H)}\right).
\end{eqnarray*}
\end{lem}

\begin{proof}
We start with the last inequality.
Define $n\in\Z$ by $2^{n} \leq \lambda < 2^{n+1}$.
If the maximum variation of $X$ in the interval $[0,t]$ is at least $\lambda$, then necessarily there exists a $k$ such that $s \leq T^{n-1}_{k-1} < T^{n-1}_k \leq s + t$.
Thus, using Lemma \ref{leftW},
\begin{align}\label{upbound}
\Pb\left(\sup_{0\leq u\leq t} |X(s+u)-X(s)|>\lambda\mid\calF_s\right)
&\leq \Pb(D_k^{n-1}<t) \nonumber\\
&= \Pb(W< \mu^{-(n-1)}t) \nonumber\\
&\leq c_2\exp\left( -c_3(t/\mu^{n-1})^{-H/(1-H)}\right)
\end{align}
Re-expressing the last inequality in terms of $\lambda$ and adjusting the constants yields the desired upper bound.

We now turn to the first inequality of the lemma.
Our proof is based on Theorem 4.3 in \cite{BarlP88}, though more work is required because our process is not Markov.
This time let $n\in\Z$ be such that $2^{n-2} \leq \lambda < 2^{n-1}$.
As before let $T^n_k(s)$ be the $k+1$-st level $n$ crossing time of the process after time $s$, and also let $T^n_{-1}(s)$ be the first level $n$ crossing time strictly before $s$.

Consider the possible level $n$ movements of the process up to time $T^n_0(s)$.
We take cases depending on the orientations of the two level $n$ crossings leading up to $T^n_0(s)$.

\begin{center}
\begin{tabular}{cccc}
(a) $--$ &
(b) $++$ &
(c) $+-$ &
(d) $-+$
\end{tabular}
\end{center}

\begin{figure}
\centering
\includegraphics[width=12cm]{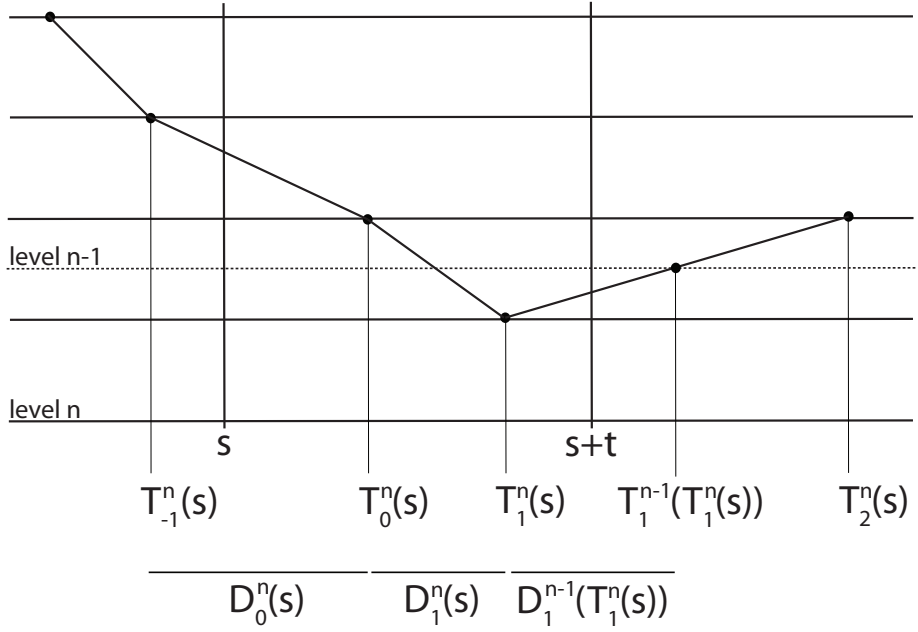}
\caption{The points considered in case (a) of the proof of Lemma \ref{cebp}.}
\label{continuity.fig}
\end{figure}

Define $\pi = \Pb(Z > 2)$ ($> 0$ by Assumption \ref{AssGW}).

In case (a) the next two level $n$ crossings will have orientation $-+$ with probability $\pi/2$.
Let $\al^n_{-1}(s)$, $\al^n_{0}(s)$, $\al^n_{1}(s)$, $\al^n_{2}(s)$, be respectively the orientations of the two crossings up to time $T^n_0(s)$ and the two crossings after time $T^n_0(s)$.
Also let $T^{n-1}_1(T^n_1(s))$ be the next level $n-1$ crossing time after $T^n_1(s)$.
Let $D^n_{0}(s) = T^n_{0}(s) - T^n_{-1}(s) \eqdist \mu^n W$, $D^n_{1}(s) = T^n_{1}(s) - T^n_{0}(s) \eqdist \mu^n W$ and $D^{n-1}_{1}(T^n_1(s)) = T^{n-1}_1(T^n_1(s)) - T^n_{1}(s) \eqdist \mu^{n-1} W$.
Note that they are independent and independent of the $\al^n_{k}(s)$.
We have
\begin{eqnarray*}
\lefteqn{ \Pb( |X(s+t)-X(s)| > \lambda \mid \calF_s, \al^n_{-1}(s) = -, \al^n_{0}(s) = -) } && \\
&\geq& \Pb( |X(s+t)-X(s)| > 2^{n-1} \mid \calF_s, \al^n_{-1}(s) = -, \al^n_{0}(s) = -) \\
&\geq& \pi/2 \Pb( |X(s+t)-X(s)| > 2^{n-1} \mid \calF_s, \\
&& \hspace{1cm} \al^n_{-1}(s) = -, \al^n_{0}(s) = -, \al^n_{1}(s) = -, \al^n_{2}(s) = +) \\
&\geq& \pi/2 \Pb( T_{1}^n(s) < s+t, T_{1}^{n-1}(T^n_1(s)) > s+t \mid \calF_s, \\
&& \hspace{1cm} \al^n_{-1}(s) = -, \al^n_{0}(s) = -, \al^n_{1}(s) = -, \al^n_{2}(s) = +) \\
&\geq& \pi/2 \Pb( T^n_{0}(s) < s + t/2, D^n_{1}(s) < t/2, D^{n-1}_{1}(T^n_1(s)) > t \mid \calF_s, \\
&& \hspace{1cm} \al^n_{-1}(s) = -, \al^n_{0}(s) = -, \al^n_{1}(s) = -, \al^n_{2}(s) = +) \\
&=& \pi/2 \Pb( T^n_{0}(s) < s + t/2 \mid \calF_s, \al^n_{-1}(s) = -, \al^n_{0}(s) = -) \\
&& \qquad \Pb( D^n_{1}(s) < t/2) \Pb( D^{n-1}_{1}(T^n_1(s)) > t)
\end{eqnarray*}
Thus from Lemmas \ref{leftW} and \ref{leftT} we have
\begin{eqnarray*}
\lefteqn{ \Pb( |X(s+t)-X(s)| > \lambda \mid \calF_s, \al^n_{-1}(s) = -, \al^n_{0}(s) = -) } && \\
&\geq&
c_1 \exp\left( -c_2(\lambda^{1/H}/t)^{H/(1-H)} \right)
\left( 1 - c_3 \exp\left( -c_4(\lambda^{1/H}/t)^{H/(1-H)} \right) \right)
\end{eqnarray*}
Choose $K$ large enough such that the last term is at least $1/2$ when $\lambda^{1/H}/t \geq K$.
Thus, since the LHS is decreasing in $\lambda$, we can find $c_5$ such that for all $t \in [0, 1]$ and $\lambda \geq 0$, the LHS is bounded below by $c_5 \exp\left( -c_2(\lambda^{1/H}/t)^{H/(1-H)} \right)$.

Cases (b) is analogous to case (a).

In case (c) we distinguish two further possibilities: (c1) the next two level $n$ crossings form an excursion (either $-+$ or $+-$); and (c2) the next two level $n$ crossings form a direct crossing (either $--$ or $++$).
In case (c1) with probability $1/2$ the excursion will be $-+$, in which case we can proceed as in case (a) to get a bound of the same form.
In case (c2) if the direct crossing is $--$ then the approach of case (a) again suffices, however if the direct crossing is $++$ then we need to modify the argument a little.
In this case we wish to bound $\Pb( |X(s+t)-X(s)| > \lambda \mid \calF_s, \al^n_{-1}(s)= +, \al^n_{0}(s) = -, \al^n_{1}(s) = +, \al^n_{2}(s) = +)$.
With probability $\pi/2$ the next pair of level $n$ crossings are the excursion $+-$.
Let $T^{n-1}_1(T^n_3(s))$ be the next level $n-1$ crossing time after $T^n_3(s)$, then we get
\begin{eqnarray*}
\lefteqn{ \Pb( |X(s+t)-X(s)| > \lambda \mid \calF_s, \al^n_{-1}(s) = +, \al^n_{0}(s) = -, \al^n_{1}(s) = +, \al^n_{2}(s) = +) } && \\
&\geq& \pi/2 \Pb( |X(s+t)-X(s)| > 2^{n-1} \mid \calF_s, \\
&&\hspace{1cm} \al^n_{-1}(s) = +, \al^n_{0,s} = -, \al^n_{1}(s) = +, \al^n_{2}(s) = +, \al^n_{3}(s) = +, \al^n_{4}(s) = -) \\
&\geq& \pi/2 \Pb( T_{3}^n(s) < s+t, T^{n-1}_1(T^n_3(s)) > s+t \mid \calF_s, \\
&&\hspace{1cm} \al^n_{-1}(s) = +, \al^n_{0,s} = -, \al^n_{1}(s) = +, \al^n_{2}(s) = +, \al^n_{3}(s) = +, \al^n_{4}(s) = -) \\
&\geq& \pi/2 \Pb( T^n_{0}(s) < s + t/4, D^n_{1}(s) < t/4, D^n_{2}(s) < t/4, D^n_{3}(s) < t/4, D^{n-1}_1(T^n_3(s)) > t \mid \calF_s, \\
&&\hspace{1cm} \al^n_{-1}(s) = +, \al^n_{0,s} = -, \al^n_{1}(s) = +, \al^n_{2}(s) = +, \al^n_{3}(s) = +, \al^n_{4}(s) = -) \\
&=& \pi/2 \Pb(T^n_{0}(s) < s + t/2 \mid \calF_s, \al^n_{-1}(s) = +, \al^n_{0}(s) = -) \\
&&\qquad \Pb( D^n_{1}(s) < t/4) \Pb( D^n_{2}(s) < t/4) \Pb( D^n_{3}(s) < t/4) \Pb(D^{n-1}_{1}(T^n_3(s)) > t).
\end{eqnarray*}
This can be bounded below in the same way as in case (a).

Case (d) is analogous to case (c).

Finally, for general $t \geq 0$, let $m$ be such that $\mu^{-m}t \leq 1$.
Then, noting that $(2^{-m}\lambda)^{1/H}/(t\mu^{-m})=\lambda^{1/H}/t\geq K$,
by the discrete scaling of $X$,
\begin{eqnarray*}
\Pb( |X(s+t)-X(s)| > \lambda \mid \calF_s)
&=& \Pb( |X(\mu^{-m}(s+t)) - X(\mu^{-m}s)| > 2^{-m}\lambda \mid \calF_{\mu^{-m}s}) \\
&\geq& c_5 \exp\left( -c_6((2^{-m}\lambda)^{1/H}/(\mu^{-m}t))^{H/(1-H)} \right) \\
&=& c_5 \exp\left( -c_6(\lambda^{1/H}/t)^{H/(1-H)} \right),
\end{eqnarray*}
which concludes the proof of the lemma.
\end{proof}

\begin{thm}\label{modcontCEBP}
Suppose that Assumptions \ref{AssGW}, \ref{AssSym} and \ref{AssZ} hold.
Let $h_H(\delta) = \delta^H|\log{\delta}|^{1-H}$, then there exist constants $c_{10}$, $c_{11} > 0$ such that
\begin{eqnarray*}
c_{10}
&\leq& \liminf_{\delta\to 0} \sup_{s,t \in [0,1], |t-s| < \delta} \frac{|X(t)-X(s)|}{h_H(t-s)} \\
&\leq& \limsup_{\delta\to 0} \sup_{s,t \in [0,1], |t-s| < \delta} \frac{|X(t)-X(s)|}{h_H(t-s)}
\ \leq\  c_{11}.
\end{eqnarray*}
\end{thm}

\begin{proof}
Consider first the lower bound.
Fix $c_1 > 0$, then for any $l > 0$ and $m = 0, 1, \ldots, 2^l-1$, put
\[
A_{m,l} = \left\{ |X((m+1)2^{-l})-X(m2^{-l})| > c_1 l^{1-H}2^{-lH} \right\}.
\]
By Lemma~\ref{cebp} we have $\Pb(A_{m,l}\mid\calF_{m2^{-l}}) \geq c_2 e^{-c_3 l}$, where $c_3 \propto c_1^{1/(1-H)}$.
By repeatedly conditioning we have
\begin{eqnarray*}
\Pb\left(\bigcap\limits_{m=0}^{2^l-1} A_{m,l}^c\right) &=& \prod_{m=0}^{2^l-1} \Pb(A^c_{m,l}\mid\calF_{m2^{-l}})
\ =\ \prod_{m=0}^{2^l-1} \left(1-\Pb(A_{m,l}\mid\calF_{m2^{-l}})\right) \\
&\leq &  \left( 1 - c_2 e^{-c_3 l}\right)^{2^l}
\ =\ \left( 1 - \frac{c_2 e^{-c_3 l}2^l}{2^l} \right)^{2^l} \\
&\leq & c_4 \exp\left( -c_2 e^{(\log 2 - c_3)l} \right)
\end{eqnarray*}
We can choose $c_1$ so that $\log 2 - c_3 > 0$, in which case the RHS above tends to $0$ as $l \to\infty$, and we have
\[
\Pb\left( |X(t+2^{-l})-X(t)| \leq c_5 h_H(2^{-l}), \;\forall l>0, \;t\in [0,1-2^{-l}] \right) = 0,
\]
which establishes the lower bound.

For the upper bound we proceed in a similar manner, though we can no longer just consider points on the lattice $2^{-l}\Z$.
For $l>0$ and $m=0,\hdots,2^l-1$, let $I_{m,l} = [m2^{-l},(m+1)2^{-l})$, and define
\begin{eqnarray*}
\Phi_{m,l} &=& \sup_{t\in I_{m,l}}|X(t)-X(m2^{-l})| \\
B_{m,l} &=& \{ \Phi_{m,l} > c_1 l^{1-H} 2^{-lH} \}
\end{eqnarray*}
From our estimate in Lemma~\ref{cebp} we have $\Pb(B_{m,l} \mid \calF_{m2^{-l}}) \leq c_2 e^{-c_3 l}$, where $c_3 \propto c_1^{1/(1-H)}$.
Thus, repeatedly conditioning on $\calF_{m2^{-l}}$ for $m = 2^l-1, \ldots, 0$, we have
\begin{eqnarray*}
\Pb( B_{m,l} \mbox{ for some } 0 \leq m < 2^l )
&=& 1 - \Pb( B_{m,l}^c \mbox{ for all } 0 \leq m < 2^l ) \\
&\leq& 1 - (1 - c_2 e^{-c_3 l})^{2^l} \\
&=& 1 - \left(1 - \frac{c_2 e^{-c_3 l}2^l}{2^l} \right)^{2^l} \\
&\leq& 1 - \exp\left( -c_4 e^{-(c_3 - \log 2)l} \right)  \\
&\leq& c_4 e^{-(c_3 - \log 2)l}.
\end{eqnarray*}
Here we have chosen $c_1$ so that $c_3 - \log 2 > 0$.

Applying the Borel-Cantelli lemma, we see that there exists an $L$ such that with probability 1
\[
\Phi_{m,l} \leq c_1 l^{1-H} 2^{-lH} \mbox{ for all } l>L \mbox{ and } 0 \leq m < 2^l.
\]

Now let $s \in I_{m,l}$ and suppose that $t$ is such that $s<t$ and $|s-t|<2^{-l}$.
Then $t\in I_{m,l}\cup I_{m+1,l}$ and we have, with probability 1,
\begin{eqnarray*}
|X(t)-X(s)|
&\leq & |X(t)-X((m+1)2^{-l})| + | X((m+1)2^{-l})-X(m2^{-l})| \\
&& \qquad + |X(m2^{-l})-X(s)| \\
&\leq & 3c_1 l^{1-H} 2^{-lH}.
\end{eqnarray*}
If we take $2^{-(l+1)} \leq \delta \leq 2^{-l}$, then we have, with probability 1,
\[
\sup_{s,t\in [0,1], |s-t|<\delta} |X(s)-X(t)|
\leq c_5 l^{1-H}2^{-lH}
\leq c_6 h_H(\delta),
\]
as required.
\end{proof}

\begin{rema} 
In the special case where the CEBP reduces to a Brownian motion, the existence of the limit in Theorem \ref{modcontCEBP} follows from Levy's modulus of continuity theorem. 
\end{rema} 

\begin{cor}
Suppose that Assumptions \ref{AssGW}, \ref{AssSym} and \ref{AssZ} hold, then the CEBP is a monofractal, in that $\Pb$-a.s. the Holder exponent $h(t) = H$ for all $t\in [0,1]$.
\end{cor}

\end{document}